\documentclass[11pt]{amsart}

\usepackage{amsmath,amsthm,amsfonts,amssymb}
\usepackage[margin=1in]{geometry}
\usepackage[bookmarks]{hyperref}

\newcommand{\set}[1]{\left\{ #1 \right\}}

\newcommand{\Hom}{\operatorname{Hom}}

\newcommand{\opH}{\operatorname{H}}
\newcommand{\Hbul}{\opH^\bullet}

\newcommand{\N}{\mathbb{N}}

\newcommand{\calV}{\mathcal{V}}

\newcommand{\frakb}{\mathfrak{b}}
\newcommand{\fraku}{\mathfrak{u}}
\newcommand{\g}{\mathfrak{g}}

\newcommand{\uz}{u_\zeta}
\newcommand{\uzg}{\uz(\g)}

\newcommand{\uzup}{\uz(\fraku^+)}

\newcommand{\uzbp}{\uz(\frakb^+)}
\newcommand{\uzbap}{\uz(\frakb_\alpha^+)}

\newcommand{\Vbp}{\calV_{\uzbp}}

\numberwithin{equation}{section}

\newtheorem{theorem}{Theorem}[section]

\theoremstyle{definition}

\title[On injective modules and support varieties for the small quantum group]{Corrigendum to ``On injective modules and support varieties for the small quantum group''}

\author{Christopher M.\ Drupieski}

\address{
Department of Mathematical Sciences \\
DePaul University \\
Chicago, IL 60614, USA}
\email{cdrupies@depaul.edu}


\thanks{The author is thankful to Alexey Sevastyanov for pointing out the flaw in the proof of \cite[Theorem 5.12]{Drupieski:2011a}.}

\begin{document}

\begin{abstract}
The proof of Theorem 5.12 in \cite{Drupieski:2011a} does not make sense as written because the algebra $u_\zeta(\frakb_\alpha^+)$ need not be a Hopf subalgebra of $\uzbp$ unless $\alpha$ is a simple root. This note describes how the proof should be modified to work around this fact.
\end{abstract}

\maketitle

All notation is taken from \cite{Drupieski:2011a}. Theorem 5.12 of \cite{Drupieski:2011a} is as follows:

\setcounter{section}{5}

\setcounter{theorem}{11}

\begin{theorem}
Let $M$ be a finite-dimensional $\uzbp$-module and let $\alpha \in \Phi^+$. Then the root vector $e_\alpha \in \fraku^+$ is an element of $\Vbp(M)$ if and only if $M$ is not projective for $u_\zeta(e_\alpha)$.
\end{theorem}

The proof given in \cite{Drupieski:2011a} involves considering the restriction of $M$ to a certain subalgebra $u_\zeta(\frakb_\alpha^+)$ of $\uzbp$. However, the proof does not make sense as written because, while $u_\zeta(\frakb_\alpha^+)$ does admit the structure of a Hopf algebra, it need not be a Hopf subalgebra of $\uzbp$ unless $\alpha$ is a simple root. The purpose of this note is to give a correct proof of the theorem.

\begin{proof}[Proof of Theorem 5.12]
Let $\alpha \in \Phi^+$ be a positive root. Write $\alpha = \sum_{\beta \in \Pi} m_\beta \beta$ as a sum of simple roots, and set $K_\alpha = \prod_{\beta \in \Pi} K_\beta^{m_\beta}$. Now let $\uzbap$ be the subalgebra of $\uzbp$ generated by $K_\alpha$ and $E_\alpha$. If $\alpha$ is a simple root, then $K_\alpha$ and $E_\alpha$ are just the defining generators of $\uzg$ labeled by $\alpha$. More generally, let $\beta \in \Pi$ be a simple root of the same length as $\alpha$. Then the assignments $E_\beta \mapsto E_\alpha$ and $K_\beta \mapsto K_\alpha$ extend to an isomorphism of algebras $u_\zeta(\frakb_\beta^+) \cong \uzbap$. The algebra $u_\zeta(\frakb_\beta^+)$ is a Hopf subalgebra of $\uzbp$, so this shows that $\uzbap$ admits the structure of a Hopf algebra (though it may not be a Hopf subalgebra of $\uzbp$ unless $\alpha$ is a simple root). In particular, taking $H = \uzbap$ and $D = u_\zeta(e_\alpha)$, it follows from \cite[Lemma 3.3]{Drupieski:2011a} that a finite-dimensional $\uzbap$-module $M$ is injective (equivalently, projective) if and only if its restriction to $u_\zeta(e_\alpha)$ is injective (respectively, projective).

Now let $M$ be a finite-dimensional $\uzbp$-module. Since $\uzbp$ is a Hopf algebra, the dual space $M^*$ is naturally a $\uzbp$-module, and then so is the tensor product $V:= M \otimes M^*$. More explicitly, the natural isomorphism $M \otimes M^* \cong \Hom_k(M,M)$ is an isomorphism of $\uzbp$-modules when the $\uzbp$-module structure on $\Hom_k(M,M)$ is defined via the formula in \cite[Remark 2.2]{Drupieski:2011a}. We define the $\uzbap$-action on $V$ to be the restriction of the $\uzbp$-action.

Next, recall that each irreducible $\uzbp$-module is one-dimensional of some weight $\lambda \in X$. More precisely, the irreducible $\uzbp$-modules are indexed by the elements of the quotient group $X/\ell X$. Denote the irreducible $\uzbp$-module of weight $\lambda$ by $k_\lambda$. Then considering a composition series for the $\uzbp$-module $M^*$, it follows that $V$ admits a $\uzbp$-module filtration with sections of the form $M \otimes k_\lambda$. The subalgebra $\uzup$ of $\uzbp$ acts trivially on $k_\lambda$, so it follows from \cite[Corollary 5.5]{Drupieski:2011a} that $V$ admits a $u_\zeta(e_\alpha)$-module filtration with each section isomorphic to $M$.

Now suppose that $M$ is projective (equivalently, injective) as a $u_\zeta(e_\alpha)$-module. Since $V$ admits a $u_\zeta(e_\alpha)$-module filtration with sections isomorphic to $M$, we conclude that $V$ is injective for $u_\zeta(e_\alpha)$. Then by \cite[Corollary 5.11]{Drupieski:2011a}, there exists $r \in \N$ such that $x_\alpha^r \in J_{\uzbp}(M)$, so $e_\alpha \notin \Vbp(M)$.

In preparation for the converse, recall that the structure of the cohomology ring $\Hbul(\uzbap,k)$ and the right $\Hbul(\uzbap,k)$-module structure of the cohomology group $\Hbul(\uzbap,V)$ depend only on the algebra structure of $\uzbap$ and on the structure of $V$ as a left $\uzbap$-module. In particular, they are independent of the existence of a Hopf algebra structure on $\uzbap$. Using the fact that $\uzbap$ admits the Hopf algebra structure of $u_\zeta(\frakb_\beta^+)$, we can deduce from the calculation of $\Hbul(\uzbp,k)$ when $\frakb^+$ is a Borel subalgebra of $\mathfrak{sl}_2$ that $\Hbul(\uzbap,k)$ is a polynomial algebra generated in cohomological degree $2$; cf.\ \cite[Theorem 5.2]{Drupieski:2011a}. Let us write $\Hbul(\uzbap,k) \cong k[e_\alpha^*]$, considering $k[e_\alpha^*]$ as the algebra of polynomial functions on the subspace of $\fraku^+$ spanned by $e_\alpha$. One can show that the restriction map $\Hbul(\uzbp,k) \rightarrow \Hbul(\uzbap,k)$ induced by the inclusion $\uzbap \hookrightarrow \uzbp$ then identifies with the natural restriction map $S(\fraku^{+*}) \rightarrow k[e_\alpha^*]$, i.e., with the map that restricts functions from $\fraku^+$ to the space $k e_\alpha$.\footnote{This description of the restriction homomorphism can be verified by calculating both cohomology rings via the argument in \cite[\S2]{Ginzburg:1993} and verifying at each step in the argument that the calculations are compatible with restriction from $\uzbp$ to $\uzbap$.}

Now suppose that $e_\alpha \notin \Vbp(M)$. By \cite[Proposition 2.4(4)]{Feldvoss:2010}, the support variety $\Vbp(M)$ is a union of relative support varieties:
\[
\Vbp(M) = \bigcup_{\lambda \in X/\ell X} \Vbp(k_\lambda,M).
\]
For each $\uzbp$-module $N$ one has $\Vbp(N,M) = \Vbp(k,M \otimes N^*)$ because of a corresponding isomorphism at the level of extension groups, so $e_\alpha \notin \Vbp(k,M \otimes k_\lambda^*)$ for each $\lambda \in X$. For the rest of the proof, redefine $V$ to be the $\uzbp$-module $M \otimes k_\lambda^*$.

Let $J_\alpha(k,M \otimes k_\lambda^*)$ be radical of the annihilator ideal for the right action of the cohomology ring $\Hbul(\uzbap,k) \cong k[e_\alpha^*]$ on $\Hbul(\uzbap,M \otimes k_\lambda^*)$. Let $\calV_\alpha(k,M \otimes k_\lambda^*)$ be the (conical) subvariety of $ke_\alpha$ defined by $J_\alpha(k,M \otimes k_\lambda^*)$. Since the restriction map $\Hbul(\uzbp,k) \rightarrow \Hbul(\uzbap,k)$ is surjective, it induces by naturality a closed embedding $\calV_\alpha(k,M \otimes k_\lambda^*) \hookrightarrow \Vbp(k,M \otimes k_\lambda^*)$. Moreover, by the assertion at the end of the previous paragraph we can identify the image of $\calV_\alpha(k,M \otimes k_\lambda^*)$ with a conical subset of $\Vbp(k,M \otimes k_\lambda^*) \cap ke_\alpha$. Since $e_\alpha \notin \Vbp(k,M \otimes k_\lambda^*)$, this conical subset must be $\set{0}$. Then we must have $\calV_\alpha(k,M \otimes k_\lambda^*) = \set{0}$ as well.

Up to this point we have considered $M \otimes k_\lambda^*$ as a $\uzbap$-module by way of restriction from $\uzbp$, with $\uzbp$ acting diagonally on $M \otimes k_\lambda^*$. On the other hand, since $\uzbap$ admits the Hopf algebra structure of $u_\zeta(\frakb_\beta^+)$, we could use the $\uzbap$-module structure on $M$, obtained via restriction from $\uzbp$, and the $\uzbap$-module structure on $k_\lambda^*$, coming from the Hopf algebra structure on $\uzbap$, together with the Hopf algebra structure on $\uzbap$ to define the diagonal action of $\uzbap$ on $M \otimes k_\lambda^* = M \otimes k_{-\lambda}$. However, in both situations one has $M \otimes k_{-\lambda} \cong M$ as a $u_\zeta(e_\alpha)$-module, from which it follows that the two $\uzbap$-actions on $M \otimes k_\lambda^*$ are the same. So now considering $\uzbap$ as a Hopf algebra, we deduce for each $\lambda \in X$ that
\[
\set{0} = \calV_\alpha(k,M \otimes k_\lambda^*) = \calV_{\uzbap}(k,M \otimes k_\lambda^*) = \calV_{\uzbap}(k_\lambda,M),
\]
and hence by \cite[Proposition 2.4(4)]{Feldvoss:2010} that $\calV_{\uzbap}(M) = \set{0}$. Then by \cite[Proposition 2.4(1)]{Feldvoss:2010}, $M$ is projective as a $\uzbap$-module, which implies that $M$ is projective for $u_\zeta(e_\alpha)$.
\end{proof}

\makeatletter
\renewcommand*{\@biblabel}[1]{\hfill#1.}
\makeatother


\providecommand{\bysame}{\leavevmode\hbox to3em{\hrulefill}\thinspace}
\providecommand{\MR}{\relax\ifhmode\unskip\space\fi MR }
\providecommand{\MRhref}[2]{%
  \href{http://www.ams.org/mathscinet-getitem?mr=#1}{#2}
}
\providecommand{\href}[2]{#2}

\end{document}